\documentclass{amsart}
\usepackage{amssymb}
\usepackage{amsmath}
\usepackage{amscd}
\usepackage[left=3cm, right=2.5cm, top=2cm, bottom=2cm]{geometry}
\usepackage{hyperref}

\newtheorem{prop}{Proposition}[section]

\newtheorem{theorem}[prop]{Theorem}

\theoremstyle{definition}

\renewcommand{\Im}{{\mathrm {Im}}}

\newcommand{\p}{{\mathfrak p}}

\DeclareFontEncoding{OT2}{}{} 

\begin{document}

\title{3-Divisibility of 9- and 27-Regular Partitions}
\author{S. Abinash}

\address{Indian Institute of Science Education and Research Thiruvananthapuram, Maruthamala P.O., Vithura, Thiruvananthapuram-695551, Kerala, India.}
\email{sarmaabinash15@iisertvm.ac.in}

\begin{abstract}
A partition of $n$ is $l$-regular if none of its parts is divisible by $l$. Let $b_l(n)$ denote the number of $l$-regular partitions of $n$. In this paper, using the theory of Hecke eigenforms explored by J.-P. Serre, we establish exact criteria for the $3$-divisibility of $b_9(n)$ and $b_{27}(n)$.
\end{abstract}

\maketitle

\section{\Large{\textbf{Introduction}}}

A partition of a positive integer $n$ is a non-increasing sequence of positive integers whose sum is $n$. The members of the sequence are called parts. For an integer $l\geq2$, a partition of $n$ is said to be $l$-regular if none of its parts is divisible by $l$. Let $b_l(n)$ denote the number of $l$-regular partitions of $n$. By convention, we assume that $b_l(0)=1$. The generating function for $b_l(n)$ is given by
\begin{align}
\sum_{n=0}^\infty b_l(n)q^n=\prod_{m=1}^\infty\frac{\left(1-q^{lm}\right)}{\left(1-q^m\right)}.\label{eq:1}
\end{align}
We also state the definition of Dedekind's eta-function,
\begin{align}
    \eta(z):=q^{1/24}\prod_{m=1}^\infty(1-q^m),\label{eq:30}
\end{align}
where $q=e^{2\pi iz}$.

In recent years, a few works have been done with regard to divisibility of $l$-regular partitions. For instance, Gordon and Ono \cite{4} showed that $b_l(n)$ is almost always divisible by $l^i$ for any prime $l$ and any positive integer $i$. Lovejoy and Penniston \cite{5} proved some results on $3$-divisibility of $b_3(n)$. For $l\in\{5,7,11 \}$, Dandurand and Penniston \cite{1} gave exact criteria for the $l$-divisibility of $b_l(n)$. In \cite{2}, Penniston established a precise description of the behaviour of $b_{11}(n)$ modulo $5$ when $5\nmid n$. 

In this paper, we establish precise criteria for $3$-divisibility of $b_9(n)$ and $b_{27}(n)$. For any positive integer $n$ and prime $p$, let $ord_p(n)$ denote the highest power of $p$ that divides $n$. The following is the first main result of this paper.

\begin{theorem}
Let $n$ be a nonnegative integer. We have that $b_9(n)$ is divisible by $3$ if and only if one of the following conditions holds.
\begin{itemize}
    \item $\#\{p:p\equiv2\pmod3, \text{ } ord_p(3n+1) \text{ is odd}\}\neq0$.
    \item $\#\{p:p\equiv1\pmod3, \text{ } ord_p(3n+1)\equiv2\pmod3\}\neq0$.
\end{itemize}
\end{theorem}

Note that for an odd prime $p$, we can write $p=x^2+y^2$ for $x,y\in\mathbb Z$ if and only if $p\equiv1\pmod 4$. Furthermore, if $p\equiv5\pmod{12}$ then we may assume, without loss of generality, that $x\equiv1\equiv-y\pmod3$ as $x^2+y^2\equiv2\pmod3$; and if $p\equiv1\pmod{12}$ then we may assume, without loss of generality, that $x\not\equiv0\pmod3$ and $y\equiv0\pmod3$ as $x^2+y^2\equiv1\pmod3$.

Suppose $p_1,\cdots,p_s$ are primes such that $p_j\equiv5\pmod{12}$ with $3\mid\mid(x_j+y_j)$ and $\alpha_j=ord_{p_j}(12n+13)$ is even. Let's denote
\begin{align}
    k&:=(-1)^{s+1}\left(\prod_{j=1}^s(-1)^{\alpha_j/2-1}p_j^{6\alpha_j-12}\right)\left(\sum_{j=1}^s\frac{\alpha_j(\alpha_j+2)}{8}\left(\prod_{i\neq j}p_i^{12}\right)\right),\label{eq:26}
\end{align}
We also denote, for each $\alpha=ord_p(12n+13)$ with $p\equiv1\pmod{12}$,
\begin{align}
    r_\alpha&:=\#\{p\equiv1\pmod{12}:ord_p(12n+13)=\alpha\},\label{eq:29}\\
    A_\alpha&:=(x_1\cdots x_{r_\alpha})^{12\alpha},\label{eq:27}\\
    B_\alpha&:=(x_1\cdots x_{r_\alpha})^{12\alpha-2}\sum_{\substack{i_1<\cdots<i_{r_\alpha-1}\\j\neq i_s}}(x_{i_1}\cdots x_{i_{r_\alpha-1}}y_j)^2.\label{eq:28}
\end{align}
Now we state the second main result of this paper.
\begin{theorem}
Let $n$ be a nonnegative integer. We have that $b_{27}(n)$ is divisible by $3$ if and only if one of the following conditions holds.
\begin{itemize}
   \item $\#\{p : p\equiv7 \text{ or } 11\pmod{12}, \text{ } ord_p(12n+13) \text{ is odd}\}\neq0$.
    \item $\#\{p : p\equiv5\pmod{12}, \text{ } ord_p(12n+13) \text{ is odd}\}\neq0,2$.
    \item $\#\{p : p\equiv5\pmod{12}, \text{ } ord_p(12n+13) \text{ is odd}\}=2;\quad 9\mid(x+y) \text{ for at least one such } p$.
    \item $\#\{p : p\equiv5\pmod{12}, \text{ } ord_p(12n+13) \text{ is odd}\}=2;\quad 3\mid\mid(x+y) \text{ for both such } p;\\
    \#\{p : p\equiv1\pmod{12}, \text{ } ord_p(12n+13)\equiv2\pmod3\}\neq0$.
    \item $\#\{p : p\equiv5\pmod{12}, \text{ } ord_p(12n+13) \text{ is odd}\}=0;\quad k\equiv0\pmod3;\\
    2\prod_\alpha(\alpha+1)^{r_\alpha}\left(\sum_\alpha\left(\alpha B_\alpha\prod_{\beta\neq\alpha}A_\beta\right)\right)\\
    \text{ }\equiv41\prod_\alpha(\alpha+1)^{r_\alpha}\left(\sum_\alpha\left(\alpha(11\alpha+26)B_\alpha\prod_{\beta\neq\alpha}A_\beta\right)\right)\pmod{3^4},\\
    \text{where $\alpha$ varies over all distinct $ord_p(12n+13)$ with $p\equiv1\pmod{12}$.}$\\
    \item $\#\{p : p\equiv5\pmod{12}, \text{ } ord_p(12n+13) \text{ is odd}\}=0;\quad k\equiv\varepsilon\pmod3 \text{ where $\varepsilon=1$ or $2$};\\
    2\prod_\alpha(\alpha+1)^{r_\alpha}\left(\sum_\alpha\left(\alpha B_\alpha\prod_{\beta\neq\alpha}A_\beta\right)\right)\\
    \text{ }\equiv41(81d+1)\prod_\alpha(\alpha+1)^{r_\alpha}\left(\sum_\alpha\left(\alpha(11\alpha+26)B_\alpha\prod_{\beta\neq\alpha}A_\beta\right)\right)\pmod{3^4},\\
    \text{where $\alpha$ varies over all distinct $ord_p(12n+13)$ with $p\equiv1\pmod{12}$ and $d=\varepsilon d'$, $d'\in\mathbb Z$ is the}\\
    \text{inverse of $\prod_{\substack{p\mid(12n+13)\\p\equiv1(12)}}(-1)^{\alpha/2}p^{6\alpha}$ in $(\mathbb Z/3^5\mathbb Z)^*$ }$.
\end{itemize}
\end{theorem}

\section{\Large{\textbf{Preliminary for Hecke Eigenform}}}

Let's consider an imaginary quadratic field $K=\mathbb Q(\sqrt d)$ with discriminant $d_K$. We denote the ring of integers of $K$ by $\mathcal O_K$. Suppose, $c$ is a Hecke character on $K$ of exponent $k-1$ and conductor $\mathfrak f_c$, where $\mathfrak f_c$ is a nonzero ideal of $\mathcal O_K$.

Now we define the series
\begin{align}
    \phi_{K,c}(z):=\sum_{\mathfrak a}c(\mathfrak a)q^{N(\mathfrak a)}\quad(q:=e^{2\pi iz}, \Im(z)>0),\label{eq:2}
\end{align}
where the sum is taken over all ideals $\mathfrak a$ of $\mathcal O_K$ which are coprime to $\mathfrak f_c$. We have that $\phi_{K,c}(z)$ is a normalized Hecke eigenform of weight $k$, level $|d_K|\cdot N(\mathfrak f_c)$ and character $\epsilon_K\omega_c$, where $\epsilon_K$ and $\omega_c$ are two Dirichlet characters defined as
\begin{align*}
    \epsilon_K(p)=\left(\frac{d}{p}\right)\text{ for primes }p\text{ such that }p\nmid2d
\end{align*}
and
\begin{align*}
    \omega_c(n)=\frac{c(n\mathcal O_K)}{n^{k-1}}\text{ for integers }n\text{ such that the ideal }n\mathcal O_K\text{ is coprime to }\mathfrak f_c.
\end{align*}
For more details please refer to \cite{1,3}.

\section{\Large{\textbf{Proof of Theorem 1.1}}}

In (\ref{eq:1}), taking $l=9$, we get
\begin{align*}
    \sum_{n=0}^\infty b_9(n)q^n=\prod_{m=1}^\infty\frac{\left(1-q^{9m}\right)}{\left(1-q^m\right)}\equiv\prod_{m=1}^\infty\left(1-q^m\right)^8\pmod3.
\end{align*}
Thus, in view of (\ref{eq:30}), we have
\begin{align}
    \sum_{n=0}^\infty b_9(n)q^{3n+1}\equiv q\cdot\prod_{m=1}^\infty\left(1-q^{3m}\right)^8\equiv\eta^8(3z)\pmod3,\label{eq:3}
\end{align}

In \cite{3}, Serre showed that $\eta^8(3z)$ is a Hecke eigenform of weight $4$ and level $9$. Let's consider $K=\mathbb Q(\sqrt{-3})$ whose discriminant is $-3$. We have $\mathcal O_K=\mathbb Z[\zeta_3]$ where $\zeta_3$ is a primitive cube root of unity. Now we fix the ideal $\mathfrak f=\sqrt{-3}\cdot\mathcal O_K$; note that $N(\mathfrak f)=3$. We define a Hecke character $c$ on $K$ of exponent $3$ and conductor $\mathfrak f$.

Suppose $\mathfrak a$ is an arbitrary ideal of $\mathcal O_K$ coprime to $\mathfrak f$. Since $\mathcal O_K$ is a PID, therefore $\mathfrak a$ is a principal ideal. Furthermore, we can always choose a generator $a$ of $\mathfrak a$ such that $a\equiv1\pmod{\mathfrak f}$. We define $c(\mathfrak a)=a^3$. It can be easily verified that $c$ thus defined is a Hecke character on $K$ of exponent $3$ and conductor $\mathfrak f$. Also, $\epsilon_K\omega_c$ turns out to be the principal character here. So, $\phi_{K,c}(z)$ from (\ref{eq:2}) is a Hecke eigenform of weight $4$ and level $9$.

Since $\eta^8(3z)=\phi_{K,c}(z)$ (see \cite{3}), therefore $\eta^8(3z)$ is a Hecke eigenform of weight $4$ and level $9$. We write
\begin{align*}
    \phi_{K,c}(z)=\sum_{n=1}^\infty s(n)q^n.
\end{align*}
By (\ref{eq:3}), we get
\begin{align*}
    \sum_{n=0}^\infty b_9(n)q^{3n+1}\equiv\sum_{n=1}^\infty s(n)q^n\pmod3.
\end{align*}
Comparing the coefficients of the terms involving $q^{3n+1}$, we get $b_9(n)\equiv s(3n+1)\pmod3$ for every nonnegative integer $n$.

Let $3n+1=\prod_{i=1}^rp_i^{e_i}$ for some $r\geq0$ where $p_i$'s are distinct primes. In view of the properties of Hecke eigenforms, we have
\begin{align*}
    s(3n+1)=\prod_{i=1}^rs(p_i^{e_i}).
\end{align*}
So, we shift our focus to $s(p^j)$ for primes $p$.

For an odd prime $p\neq3$, note that
\begin{align*}
    p\text{ is totally split in }K&\Longleftrightarrow\left(\frac{-3}{p}\right)=1\\
    &\Longleftrightarrow p\text{ is represented by some primitive form of discriminant }-12\\
    &\Longleftrightarrow p=x^2+3y^2\text{ for some }x,y\in\mathbb Z\\
    &\Longleftrightarrow p\equiv1\pmod3.
\end{align*}
Consequently, $p$ is inert in $K$ $\Longleftrightarrow$ $p\equiv2\pmod3$. Also, $2$ is inert in $K$.

Let $p\equiv2\pmod3$. Since $p$ is inert, therefore no ideal of $\mathcal O_K$ has norm $p^j$ if $j$ is odd. Thus, $s(p^j)=0$ if $j$ is odd. If $j$ is even, then $p^{j/2}\mathcal O_K$ is the only ideal with norm $p^j$ and by definition, clearly $s(p^j)\not\equiv0\pmod3$ in this case.

Let $p\equiv1\pmod3$. We can write $p=x^2+3y^2$ for some $x,y\in\mathbb Z$ with $x\equiv1\pmod3$. Only two ideals of $\mathcal O_K$ have norm $p$, namely, $(x+y\sqrt{-3})\mathcal O_K$ and $(x-y\sqrt{-3})\mathcal O_K$. So,
\begin{align*}
    s(p)&=c\left((x+y\sqrt{-3})\mathcal O_K\right)+c\left((x-y\sqrt{-3})\mathcal O_K\right)\\
    &=(x+y\sqrt{-3})^3+(x-y\sqrt{-3})^3 \quad\quad [x\equiv1\pmod3\Longrightarrow x\pm y\sqrt{-3}\equiv1\pmod{\mathfrak f}]\\
    &=4x^3-12y^2x-px\\
    &\equiv2\pmod3.
\end{align*}
Also, $s(1)=1\equiv1\pmod3$. Now, by the properties of Hecke eigenforms, we have
\begin{align*}
    s(p^2)&=s(p)s(p)-p^3s(1)\\
    &\equiv2\cdot2-1\cdot1\pmod3\\
    &\equiv0\pmod3.
\end{align*}
Using mathematical induction, it can be easily proved that for any $k\geq0$,
\begin{align*}
    s(p^{3k})\equiv1\pmod3,\quad s(p^{3k+1})\equiv2\pmod3,\quad s(p^{3k+2})\equiv0\pmod3.
\end{align*}
Hence, the proof is completed.

\section{\Large{\textbf{Proof of Theorem 1.2}}}

In (\ref{eq:1}), taking $l=27$, we get
\begin{align*}
    \sum_{n=0}^\infty b_{27}(n)q^n=\prod_{m=1}^\infty\frac{\left(1-q^{27m}\right)}{\left(1-q^m\right)}\equiv\prod_{m=1}^\infty\left(1-q^m\right)^{26}\pmod3.
\end{align*}
Thus, in view of (\ref{eq:30}), we have
\begin{align}
    \sum_{n=0}^\infty b_{27}(n)q^{12n+13}\equiv q^{13}\cdot\prod_{m=1}^\infty\left(1-q^{12m}\right)^{26}\equiv\eta^{26}(12z)\pmod3.\label{eq:4}
\end{align}

In \cite{3}, Serre showed that $\eta^{26}(12z)$ is a linear combination of Hecke eigenforms of weight $13$, level $144$ and character $\chi_{-1}(\bullet):=\left(\frac{-1}{\bullet}\right)$.

Firstly, let's consider $K_1=\mathbb Q(\sqrt{-3})$ whose discriminant is $-3$. We have $\mathcal O_1=\mathcal O_{K_1}=\mathbb Z[\zeta_3]$ where $\zeta_3$ is a primitive cube root of unity. Now we fix the ideal $\mathfrak f_1=4\sqrt{-3}\cdot\mathcal O_1$; note that $N(\mathfrak f_1)=48$. We define Hecke characters $c_{1\pm}$ on $K_1$ of exponent $12$ and conductor $\mathfrak f_1$.

Suppose $\mathfrak a$ is an arbitrary ideal of $\mathcal O_1$ coprime to $\mathfrak f_1$. Since $\mathcal O_1$ is a PID, therefore $\mathfrak a$ is a principal ideal. Furthermore, we can always choose a generator $a=x+y\sqrt{-3}$ of $\mathfrak a$ with $x,y\in\mathbb Z$, $x+y\equiv1\pmod2$ and $x\equiv1\pmod3$. We define
\begin{align*}
    c_{1\pm}(\mathfrak a)=(-1)^{(x\mp y-1)/2}a^{12}.
\end{align*}
It can be easily verified that $c_{1\pm}$ thus defined are Hecke characters on $K_1$ of exponent $12$ and conductor $\mathfrak f_1$. Also, $\epsilon_{K_1}\omega_{c_{1+}}=\chi_{-1}=\epsilon_{K_1}\omega_{c_{1-}}$. So, $\phi_{K_1,c_{1\pm}}(z)$ from (\ref{eq:2}) are Hecke eigenforms of weight $13$, level $144$ and character $\chi_{-1}$.

Next, let's consider $K_2=\mathbb Q(i)$ whose discriminant is $-4$. We have $\mathcal O_2=\mathcal  O_{K_2}=Z[i]$. Now we fix the ideal $\mathfrak f_2=6\cdot\mathcal O_2$; note that $N(\mathfrak f_2)=36$. We define Hecke characters $c_{2\pm}$ on $K_2$ of exponent $12$ and conductor $\mathfrak f_2$.

Suppose $\mathfrak b$ is an arbitrary ideal of $\mathcal O_2$ coprime to $\mathfrak f_2$. Since $\mathcal O_2$ is a PID, therefore $\mathfrak b$ is a principal ideal. For a generator $b$ of $\mathfrak b$, we determine $m\in\mathbb Z/2\mathbb Z$ and $n\in\mathbb Z/8\mathbb Z$ such that
\begin{align*}
    b\equiv i^m\pmod{2\cdot\mathcal O_2} \quad \text{and} \quad b\equiv(1-i)^n\pmod{3\cdot\mathcal O_2}.
\end{align*}
We define
\begin{align*}
    c_{2\pm}(\mathfrak b)=(-1)^{3m}(\pm i)^{3n}b^{12}.
\end{align*}
It can be easily verified that $c_{2\pm}$ thus defined are Hecke characters on $K_2$ of exponent $12$ and conductor $\mathfrak f_2$. Also, $\epsilon_{K_2}\omega_{c_{2+}}=\chi_{-1}=\epsilon_{K_2}\omega_{c_{2-}}$. So, $\phi_{K_2,c_{2\pm}}(z)$ from (\ref{eq:2}) are Hecke eigenforms of weight $13$, level $144$ and character $\chi_{-1}$.

Since
\begin{align}
    \eta^{26}(12z)=\frac{1}{32617728}\left(\phi_{K_1,c_{1+}}(z)+\phi_{K_1,c_{1-}}(z)-\phi_{K_2,c_{2+}}(z)-\phi_{K_2,c_{2-}}(z)\right) \quad \text{ (see \cite{3})},\label{eq:9}
\end{align}
therefore $\eta^{26}(12z)$ is a linear combination of Hecke eigenforms of weight $13$, level $144$ and character $\chi_{-1}$. We write
\begin{align*}
    \phi_{K_1,c_{1\pm}}(z)=\sum_{n=1}^\infty t_{1\pm}(n)q^n, \quad \phi_{K_2,c_{2\pm}}(z)=\sum_{n=1}^\infty t_{2\pm}(n)q^n.
\end{align*}
By (\ref{eq:4}), we get
\begin{align*}
    \sum_{n=0}^\infty b_{27}(n)q^{12n+13}\equiv\frac{1}{32617728}\left(\sum_{n=1}^\infty\big(t_{1+}(n)+t_{1-}(n)-t_{2+}(n)-t_{2-}(n)\big)q^n\right)\pmod3.
\end{align*}
Comparing the coefficients of the terms involving $q^{12n+13}$, we get
\begin{align*}
b_{27}(n)\equiv\frac{1}{32617728}\big(t_{1+}(12n+13)+t_{1-}(12n+13)-t_{2+}(12n+13)-t_{2-}(12n+13)\big)\pmod3
\end{align*}
for every nonnegative integer $n$. Since $3^4\mid\mid32617728$, therefore
\begin{align*}
b_{27}(n)\equiv0\pmod3 \Longleftrightarrow t_{1+}(12n+13)+t_{1-}(12n+13)\equiv t_{2+}(12n+13)+t_{2-}(12n+13)\pmod{3^5}.
\end{align*}
In view of the properties of Hecke eigenforms, we have
\begin{align*}
    t_{1\pm}(12n+13)=\prod_{p\mid(12n+13)}t_{1\pm}(p^\alpha),\quad t_{2\pm}(12n+13)=\prod_{p\mid(12n+13)}t_{2\pm}(p^\alpha).
\end{align*}

Before proceeding further, let's have a look at $t_{1+}(p^\alpha)$, $t_{1-}(p^\alpha)$, $t_{2+}(p^\alpha)$, $t_{2-}(p^\alpha)$ for a prime $p\neq3$. Let $t$ denote any one among $t_{1\pm}$ and $t_{2\pm}$. In view of the properties of Hecke eigenforms, if $p\equiv1\pmod4$ then we have for any $\beta\geq0$,
\begin{align}
    t(p^{2\beta})&=t(p)^{2\beta}-\binom{2\beta-1}{1}p^{12}t(p)^{2\beta-2}+\cdots+(-1)^{\beta-1}\binom{\beta+1}{2}p^{12\beta-12}t(p)^2+(-1)^\beta p^{12\beta},\label{eq:5}\\
    t(p^{2\beta+1})&=t(p)^{2\beta+1}-\binom{2\beta}{1}p^{12}t(p)^{2\beta-1}+\cdots+(-1)^{\beta-1}\binom{\beta+2}{3}p^{12\beta-12}t(p)^3+(-1)^\beta\binom{\beta+1}{1}p^{12\beta}t(p).\label{eq:6}
\end{align}
If $p\equiv3\pmod4$ then we have for any $\beta\geq0$,
\begin{align}
    t(p^{2\beta})&=t(p)^{2\beta}+\binom{2\beta-1}{1}p^{12}t(p)^{2\beta-2}+\cdots+\binom{\beta+1}{2}p^{12\beta-12}t(p)^2+p^{12\beta},\label{eq:7}\\
    t(p^{2\beta+1})&=t(p)^{2\beta+1}+\binom{2\beta}{1}p^{12}t(p)^{2\beta-1}+\cdots+\binom{\beta+2}{3}p^{12\beta-12}t(p)^3+\binom{\beta+1}{1}p^{12\beta}t(p).\label{eq:8}
\end{align}

Firstly, let $p\equiv11\pmod{12}$. Since $p\equiv2\pmod3$ and $p\equiv3\pmod4$, therefore $p$ is inert in $K_1$ and $K_2$, respectively. Consequently, $t_{1\pm}(p)=0=t_{2\pm}(p)$. By (\ref{eq:7}) and (\ref{eq:8}), we have
\begin{align}\label{eq:10}
    t_{1\pm}(p^\alpha)=t_{2\pm}(p^\alpha)=
    \begin{cases}
    p^{6\alpha} &\text{ if $\alpha$ is even,}\\
    0 &\text{ if $\alpha$ is odd.}
    \end{cases}
\end{align}

Secondly, let $p\equiv7\pmod{12}$. Since $p\equiv3\pmod4$, therefore $p$ is inert in $K_2$. Consequently, $t_{2\pm}(p)=0$. By (\ref{eq:7}) and (\ref{eq:8}), we have
\begin{align}\label{eq:11}
    t_{2\pm}(p^\alpha)=
    \begin{cases}
    p^{6\alpha} &\text{ if $\alpha$ is even,}\\
    0 &\text{ if $\alpha$ is odd.}
    \end{cases}
\end{align}
Since $p\equiv1\pmod3$, therefore we can write $p=z^2+3w^2=(z+w\sqrt{-3})(z-w\sqrt{-3})$ for some $z,w\in\mathbb Z$ with $z\equiv1\pmod4$. Moreover, as $z^2+3w^2\equiv3\pmod4$, we must have $z$ even and $w$ odd and thus we can choose $z,w$ such that $z+w\equiv1\pmod4$ and $z-w\equiv3\pmod4$. Hence, we have
\begin{align}\label{eq:15}
    t_{1+}(p)=-(z+w\sqrt{-3})^{12}+(z-w\sqrt{-3})^{12}=-t_{1-}(p).
\end{align}
By (\ref{eq:8}), for any odd $\alpha$, $t_{1+}(p^\alpha)=-t_{1-}(p^\alpha)$. Further, it can be easily verified that
\begin{align}\label{eq:14}
    t_{1+}(p)^2=t_{1-}(p)^2\equiv0\pmod{3^5}.
\end{align}
By (\ref{eq:7}), we have for any even $\alpha$
\begin{align}\label{eq:16}
    t_{1\pm}(p^\alpha)\equiv p^{6\alpha}\pmod{3^5},
\end{align}

Thirdly, let $p\equiv5\pmod{12}$. Since $p\equiv2\pmod3$, therefore $p$ is inert in $K_1$.Consequently, $t_{1\pm}(p)=0$. By (\ref{eq:5}) and (\ref{eq:6}), we have
\begin{align}\label{eq:12}
    t_{1\pm}(p^\alpha)=
    \begin{cases}
    (-1)^{\alpha/2}p^{6\alpha} &\text{ if $\alpha$ is even,}\\
    0 &\text{ if $\alpha$ is odd.}
    \end{cases}
\end{align}
Since $p\equiv1\pmod4$, therefore we can write $p=x^2+y^2=(x+iy)(x-iy)$ for some $x,y\in\mathbb Z$ with $x$ odd and $y$ even. Moreover, as $x^2+y^2\equiv2\pmod3$, we must have $x\not\equiv0\pmod3$ and $y\not\equiv0\pmod3$ and thus we can choose $x,y$ such that $x\equiv1\equiv-y\pmod3$. Hence, we have
\begin{align}\label{eq:19}
    t_{2+}(p)=-i(x+iy)^{12}+i(x-iy)^{12}=-t_{2-}(p).
\end{align}
By (\ref{eq:6}), for any odd $\alpha$, $t_{2+}(p^\alpha)=-t_{2-}(p^\alpha)$. Further, it can be easily verified that
\begin{align}\label{eq:18}
    t_{2+}(p)^2=t_{2-}(p)^2\equiv
    \begin{cases}
    162\pmod{3^5} &\text{if $3\mid\mid(x+y)$,}\\
    0\pmod{3^5} &\text{if $9\mid(x+y)$.}
    \end{cases}
\end{align}
By (\ref{eq:5}), we have for any even $\alpha$
\begin{align}\label{eq:20}
    t_{2\pm}(p^\alpha)\equiv
    \begin{cases}
    (-1)^{\alpha/2-1}\left(\left(\frac{(\alpha+2)\alpha}{8}\right)p^{6\alpha-12}162-p^{6\alpha}\right)\pmod{3^5} &\text{if $3\mid\mid(x+y)$,}\\
    (-1)^{\alpha/2}p^{6\alpha}\pmod{3^5} &\text{if $9\mid(x+y)$,}
    \end{cases}
\end{align}

Finally, let $p\equiv1\pmod{12}$. Since $p\equiv1\pmod3$, therefore we can write $p=z^2+3w^2=(z+w\sqrt{-3})(z-w\sqrt{-3})$ for some $z.w\in\mathbb Z$ with $z\equiv1\pmod3$. Moreover, as $z^2+3w^2\equiv1\pmod4$, we must have $z$ odd and $w$ even and thus we have $z+w\equiv z-w\pmod4$. Hence, we have
\begin{align}\label{eq:13}
    t_{1\pm}(p)=
    \begin{cases}
    (z+w\sqrt{-3})^{12}+(z-w\sqrt{-3})^{12} &\text{if $z+w\equiv1\equiv z-w\pmod4$,}\\
    -(z+w\sqrt{-3})^{12}-(z-w\sqrt{-3})^{12} &\text{if $z+w\equiv3\equiv z-w\pmod4$.}
    \end{cases}
\end{align}
Since $p\equiv1\pmod4$, therefore we can write $p=x^2+y^2=(x+iy)(x-iy)$ for some $x,y\in\mathbb Z$ with either $x$ odd, $y$ even or $x$ even, $y$ odd. Moreover, as $x^2+y^2\equiv1\pmod3$, we can choose $x\equiv1\pmod3$ and $y\equiv0\pmod3$. Hence, we have
\begin{align}\label{eq:17}
    t_{2\pm}(p)=
    \begin{cases}
    (x+iy)^{12}+(x-iy)^{12} &\text{if $x$ odd, $y$ even,}\\
    -(x+iy)^{12}-(x-iy)^{12} &\text{if $x$ even, $y$ odd.}
    \end{cases}
\end{align}
From (\ref{eq:9}), it is clear that $t_{1+}(p)+t_{1-}(p)-t_{2+}(p)-t_{2-}(p)\equiv0\pmod{3^4}$. Note that
\begin{align*}
    \pm\left((z+w\sqrt{-3})^{12}+(z-w\sqrt{-3})^{12}+(x+iy)^{12}+(x-iy)^{12}\right)\equiv\pm2(x^{12}+z^{12})\equiv\pm1\not\equiv0\pmod3.
\end{align*}
Thus, we must have if $x$ odd, $y$ even then $z+w\equiv1\equiv z-w\pmod4$ and if $x$ even, $y$ odd then $z+w\equiv3\equiv z-w\pmod4$. In view of the equation $x^2+y^2=z^2+3w^2$, we can conclude that
\begin{align*}
    t_{1\pm}(p)^\alpha\equiv((-1)^y)^\alpha\left(2^\alpha x^{12\alpha}+12\alpha2^{\alpha-1}x^{12\alpha-2}y^2\right)\pmod{3^5}.
\end{align*}
We can easily conclude that
\begin{align*}
    t_{2\pm}(p)^\alpha\equiv((-1)^y)^\alpha\left(2^\alpha x^{12\alpha}+111\alpha2^{\alpha-1}x^{12\alpha-2}y^2\right)\pmod{3^5}.
\end{align*}
Also, for any $n$
\begin{align*}
    p^{12n}=(x^2+y^2)^{12n}\equiv x^{24n}+12nx^{24n-2}y^2\pmod{3^5}.
\end{align*}
Now by (\ref{eq:5}), we have for $\alpha=2\beta$,
\begin{align}
    t_{1\pm}(p^\alpha)&\equiv\left(\sum_{\delta=0}^\beta(-1)^\delta\binom{2\beta-\delta}{\delta}2^{2(\beta-\delta)}\right)x^{24\beta}+\left(\sum_{\delta=0}^\beta(-1)^\delta\binom{2\beta-\delta}{\delta}2^{2(\beta-\delta)}\right)12\beta x^{24\beta-2}y^2\pmod{3^5}\label{eq:21}\\
    &\equiv\left(\alpha+1\right)x^{12\alpha}+6\alpha\left(\alpha+1\right)x^{12\alpha-2}y^2\pmod{3^5},\notag\\
    t_{2\pm}(p^\alpha)&\equiv\left(\sum_{\delta=0}^\beta(-1)^\delta\binom{2\beta-\delta}{\delta}2^{2(\beta-\delta)}\right)x^{24\beta}\label{eq:22}\\
    &\quad+\left(\sum_{\delta=0}^\beta(-1)^\delta\binom{2\beta-\delta}{\delta}2^{2(\beta-\delta)}(111(\beta-\delta)+12\delta)\right)x^{24\beta-2}y^2\pmod{3^5}\notag\\
    &\equiv\left(\alpha+1\right)x^{12\alpha}+123\alpha(\alpha+1)(11\alpha+26)x^{12\alpha-2}y^2\pmod{3^5}\notag.
\end{align}
And by (\ref{eq:6}), we have for $\alpha=2\beta+1$,
\begin{align}
    t_{1\pm}(p^\alpha)&\equiv(-1)^y\Bigg(\left(\sum_{\delta=0}^\beta(-1)^\delta\binom{2\beta+1-\delta}{\delta}2^{2(\beta-\delta)+1}\right)x^{24\beta+12}\label{eq:23}\\
    &\quad+\left(\sum_{\delta=0}^\beta(-1)^\delta\binom{2\beta+1-\delta}{\delta}2^{2(\beta-\delta)}\right)12(2\beta+1)x^{24\beta+10}y^2\Bigg)\pmod{3^5}\notag\\
    &\equiv(-1)^y\left((\alpha+1)x^{12\alpha}+6\alpha(\alpha+1)x^{12\alpha-2}y^2\right)\pmod{3^5},\notag\\
    t_{2\pm}(p^\alpha)&\equiv(-1)^y\Bigg(\left(\sum_{\delta=0}^\beta(-1)^\delta\binom{2\beta+1-\delta}{\delta}2^{2(\beta-\delta)+1}\right)x^{24\beta+12}\label{eq:24}\\
    &\quad+\left(\sum_{\delta=0}^\beta(-1)^\delta\binom{2\beta+1-\delta}{\delta}2^{2(\beta-\delta)}(111(2\beta+1-2\delta)+12(2\delta))\right)x^{24\beta+10}y^2\Bigg)\pmod{3^5}\notag\\
    &\equiv(-1)^y\left((\alpha+1)x^{12\alpha}+123\alpha(\alpha+1)(11\alpha+26)x^{12\alpha-2}y^2\right)\pmod{3^5}\notag.
\end{align}

We now proceed to complete the proof. Let there exist a prime $p\equiv11\pmod{12}$ with $ord_p(12n+13)$ odd. By (\ref{eq:10}), it is clear that $b_{27}(n)\equiv0\pmod3$. Assume $ord_p(12n+13)$ is even for all $p\equiv11\pmod{12}$. Again by (\ref{eq:10}), we have
\begin{align*}
    \prod_{\substack{p\mid(12n+13)\\p\equiv11(12)}}t_{1\pm}(p^\alpha)=\prod_{\substack{p\mid(12n+13)\\p\equiv11(12)}}p^{6\alpha}=\prod_{\substack{p\mid(12n+13)\\p\equiv11(12)}}t_{2\pm}(p^\alpha),
\end{align*}
and $3\nmid p^{6\alpha}$ for all $p\equiv11\pmod{12}$. Thus,
\begin{align*}
    &b_{27}(n)\equiv0\pmod3\\
    &\Longleftrightarrow\prod_{\substack{p\mid(12n+13)\\p\not\equiv11(12)}}t_{1+}(p^\alpha)+\prod_{\substack{p\mid(12n+13)\\p\not\equiv11(12)}}t_{1-}(p^\alpha)\equiv\prod_{\substack{p\mid(12n+13)\\p\not\equiv11(12)}}t_{2+}(p^\alpha)+\prod_{\substack{p\mid(12n+13)\\p\not\equiv11(12)}}t_{2-}(p^\alpha)\pmod{3^5}.
\end{align*}

Let there exist a prime $p\equiv7\pmod{12}$ with $ord_p(12n+13)$ odd. By (\ref{eq:11}), it is clear that
\begin{align*}
    b_{27}(n)\equiv0\pmod3\Longleftrightarrow\prod_{\substack{p\mid(12n+13)\\p\not\equiv11(12)}}t_{1+}(p^\alpha)+\prod_{\substack{p\mid(12n+13)\\p\not\equiv11(12)}}t_{1-}(p^\alpha)\equiv0\pmod{3^5}.
\end{align*}
Note, in view of (\ref{eq:12}) and (\ref{eq:13}), that if $p\equiv1,5\pmod{12}$ then $t_{1+}(p^\alpha)=t_{1-}(p^\alpha)$ for any $\alpha$. Thus, we have
\begin{align*}
    \prod_{\substack{p\mid(12n+13)\\p\not\equiv11(12)}}t_{1+}(p^\alpha)+\prod_{\substack{p\mid(12n+13)\\p\not\equiv11(12)}}t_{1-}(p^\alpha)=\prod_{\substack{p\mid(12n+13)\\p\equiv1,5(12)}}t_{1+}(p^\alpha)\left(\prod_{\substack{p\mid(12n+13)\\p\equiv7(12)}}t_{1+}(p^\alpha)+\prod_{\substack{p\mid(12n+13)\\p\equiv7(12)}}t_{1-}(p^\alpha)\right).
\end{align*}
Also note, in view of (\ref{eq:14}), that if $p\equiv7\pmod{12}$ then $t_{1+}(p^\alpha)=t_{1-}(p^\alpha)$ for any even $\alpha$. As $t_{1+}(p^\alpha)=-t_{1-}(p^\alpha)$ for any odd $\alpha$, we have
\begin{align*}
    \prod_{\substack{p\mid(12n+13)\\p\equiv7(12)}}t_{1+}(p^\alpha)+\prod_{\substack{p\mid(12n+13)\\p\equiv7(12)}}t_{1-}(p^\alpha)&=\prod_{\substack{p\mid(12n+13)\\p\equiv7(12)\\\alpha \text{ even}}}t_{1+}(p^\alpha)\left(\prod_{\substack{p\mid(12n+13)\\p\equiv7(12)\\\alpha \text{ odd}}}t_{1+}(p^\alpha)+\prod_{\substack{p\mid(12n+13)\\p\equiv7(12)\\\alpha \text{ odd}}}t_{1-}(p^\alpha)\right)\\
    &=(1+(-1)^u)\prod_{\substack{p\mid(12n+13)\\p\equiv7(12)\\\alpha \text{ even}}}t_{1+}(p^\alpha)\prod_{\substack{p\mid(12n+13)\\p\equiv7(12)\\\alpha \text{ odd}}}t_{1+}(p^\alpha),
\end{align*}
where $u$ is the number of odd $\alpha$. If $u$ is odd, clearly $b_{27}(n)\equiv0\pmod3$. Suppose $u$ is even. Using (\ref{eq:15}) it can be verified that for any two primes $p_1,p_2\equiv7\pmod{12}$, $t_{1+}(p_1)t_{1+}(p_2)\equiv0\pmod{3^5}$. Thus, we get that $b_{27}(n)\equiv0\pmod{3}$. Assume $ord_p(12n+13)$ is even for all $p\equiv7\pmod{12}$. By (\ref{eq:11}) and (\ref{eq:16}), we have
\begin{align*}
    \prod_{\substack{p\mid(12n+13)\\p\equiv7(12)}}t_{1\pm}(p^\alpha)\equiv\prod_{\substack{p\mid(12n+13)\\p\equiv7(12)}}p^{6\alpha}\equiv\prod_{\substack{p\mid(12n+13)\\p\equiv7(12)}}t_{2\pm}(p^\alpha)\pmod{3^5},
\end{align*}
and $3\nmid p^{6\alpha}$ for all $p\equiv7\pmod{12}$. Thus,
\begin{align*}
    &b_{27}(n)\equiv0\pmod3\\
    &\Longleftrightarrow\prod_{\substack{p\mid(12n+13)\\p\equiv1,5(12)}}t_{1+}(p^\alpha)+\prod_{\substack{p\mid(12n+13)\\p\equiv1,5(12)}}t_{1-}(p^\alpha)\equiv\prod_{\substack{p\mid(12n+13)\\p\equiv1,5(12)}}t_{2+}(p^\alpha)+\prod_{\substack{p\mid(12n+13)\\p\equiv1,5(12)}}t_{2-}(p^\alpha)\pmod{3^5}.
\end{align*}

Let there exist a prime $p\equiv5\pmod{12}$ with $ord_p(12n+13)$ odd. By (\ref{eq:12}), it is clear that
\begin{align*}
    b_{27}(n)\equiv0\pmod3\Longleftrightarrow\prod_{\substack{p\mid(12n+13)\\p\equiv1,5(12)}}t_{2+}(p^\alpha)+\prod_{\substack{p\mid(12n+13)\\p\equiv1,5(12)}}t_{2-}(p^\alpha)\equiv0\pmod{3^5}.
\end{align*}
Note, in view of (\ref{eq:17}), that if $p\equiv1\pmod{12}$ then $t_{2+}(p^\alpha)=t_{2-}(p^\alpha)$ for any $\alpha$. Thus we have
\begin{align*}
    \prod_{\substack{p\mid(12n+13)\\p\equiv1,5(12)}}t_{2+}(p^\alpha)+\prod_{\substack{p\mid(12n+13)\\p\equiv1,5(12)}}t_{2-}(p^\alpha)=\prod_{\substack{p\mid(12n+13)\\p\equiv1(12)}}t_{2+}(p^\alpha)\left(\prod_{\substack{p\mid(12n+13)\\p\equiv5(12)}}t_{2+}(p^\alpha)+\prod_{\substack{p\mid(12n+13)\\p\equiv5(12)}}t_{2-}(p^\alpha)\right).
\end{align*}
Also note, in view of (\ref{eq:18}), that if $p\equiv5\pmod{12}$ then $t_{2+}(p^\alpha)=t_{2-}(p^\alpha)$ for any even $\alpha$. As $t_{1+}(p^\alpha)=-t_{1-}(p^\alpha)$ for any odd $\alpha$, we have
\begin{align*}
    \prod_{\substack{p\mid(12n+13)\\p\equiv5(12)}}t_{2+}(p^\alpha)+\prod_{\substack{p\mid(12n+13)\\p\equiv5(12)}}t_{2-}(p^\alpha)&=\prod_{\substack{p\mid(12n+13)\\p\equiv5(12)\\\alpha \text{ even}}}t_{2+}(p^\alpha)\left(\prod_{\substack{p\mid(12n+13)\\p\equiv5(12)\\\alpha \text{ odd}}}t_{2+}(p^\alpha)+\prod_{\substack{p\mid(12n+13)\\p\equiv5(12)\\\alpha \text{ odd}}}t_{2-}(p^\alpha)\right)\\
    &=(1+(-1)^u)\prod_{\substack{p\mid(12n+13)\\p\equiv5(12)\\\alpha \text{ even}}}t_{2+}(p^\alpha)\prod_{\substack{p\mid(12n+13)\\p\equiv5(12)\\\alpha \text{ odd}}}t_{2+}(p^\alpha),
\end{align*}
where $u$ is the number of odd $\alpha$. If $u$ is odd, clearly $b_{27}(n)\equiv0\pmod3$. Suppose $u$ is even. Using (\ref{eq:19}) it can be verified that for any two primes $p_1,p_2\equiv5\pmod{12}$, when written as $p_i=x_i^2+y_i^2$ for $i=1,2$,
\begin{align*}
    t_{2+}(p_1)t_{2+}(p_2)\equiv
    \begin{cases}
    81 \text{ or } 162 &\text{if $3\mid\mid(x_i+y_1)$ for both $i=1,2$,}\\
    0 &\text{if either $9\mid(x_1+y_1)$ or $9\mid(x_2+y_2)$.}
    \end{cases}
\end{align*}
Thus, if $u\geq4$ then $b_{27}(n)\equiv0\pmod3$. For $u=2$, if either $9\mid(x_1+y_1)$ or $9\mid(x_2+y_2)$ then also $b_{27}(n)\equiv0\pmod3$. Suppose $u=2$ and $3\mid\mid(x_i+y_1)$ for both $i=1,2$. In view of (\ref{eq:20}), we get
\begin{align*}
    \prod_{\substack{p\mid(12n+13)\\p\equiv5(12)\\\alpha \text{ even}}}t_{2+}(p^\alpha)\not\equiv0\pmod3.
\end{align*}
Hence, in view of (\ref{eq:22}) and (\ref{eq:24}), we get
\begin{align*}
  b_{27}(n)\equiv0\pmod3&\Longleftrightarrow\prod_{\substack{p\mid(12n+13)\\p\equiv1(12)}}t_{2+}(p^\alpha)\equiv0\pmod3\\
  &\Longleftrightarrow\exists\quad p\equiv1\pmod{12} \quad \text{with} \quad ord_p(12n+13)\equiv2\pmod3.
\end{align*}
Assume $ord_p(12n+13$ is even for all $p\equiv5\pmod{12}$. By (\ref{eq:12}) and (\ref{eq:20}), we have
\begin{align}
    \prod_{\substack{p\mid(12n+13)\\p\equiv5(12)}}t_{1\pm}(p^\alpha)&=\prod_{\substack{p\mid(12n+13)\\p\equiv5(12)}}(-1)^{\alpha/2}p^{6\alpha},\notag\\
    \prod_{\substack{p\mid(12n+13)\\p\equiv5(12)}}t_{2\pm}(p^\alpha)&\equiv\underbrace{\prod_{\substack{p\mid(12n+13)\\p\equiv5(12)}}\left((-1)^{\alpha/2-1}\frac{(\alpha+2)\alpha}{8}p^{6\alpha-12}162+(-1)^{\alpha/2}p^{6\alpha}\right)}_\text{those $p$ for which $3\mid\mid(x+y)$}\label{eq:25}\\
    &\quad\times\underbrace{\prod_{\substack{p\mid(12n+13)\\p\equiv5(12)}}(-1)^{\alpha/2}p^{6\alpha}}_\text{those $p$ for which $9\mid(x+y)$}\pmod{3^5}\notag\\
    &\equiv k\cdot81+\prod_{\substack{p\mid(12n+13)\\p\equiv5(12)}}(-1)^{\alpha/2}p^{6\alpha}\pmod{3^5}\notag,
\end{align}
where $k$ is defined in (\ref{eq:26}).

Assume $k\equiv0\pmod3$. We have
\begin{align*}
    &b_{27}(n)\equiv0\pmod3\\
    \Longleftrightarrow&\prod_{\substack{p\mid(12n+13)\\p\equiv1(12)}}t_{1+}(p^\alpha)\prod_{\substack{p\mid(12n+13)\\p\equiv5(12)}}(-1)^{\alpha/2}p^{6\alpha}\equiv\prod_{\substack{p\mid(12n+13)\\p\equiv1(12)}}t_{2+}(p^\alpha)\prod_{\substack{p\mid(12n+13)\\p\equiv5(12)}}(-1)^{\alpha/2}p^{6\alpha}\pmod{3^5}\\
    \Longleftrightarrow&\prod_{\substack{p\mid(12n+13)\\p\equiv1(12)}}t_{1+}(p^\alpha)\equiv\prod_{\substack{p\mid(12n+13)\\p\equiv1(12)}}t_{2+}(p^\alpha)\pmod{3^5}.
\end{align*}
Suppose $ord_{p_1}(12n+13)=ord_{p_2}(12n+13)=\cdots=ord_{p_r}(12n+13)=\alpha$, $p_l\equiv1\pmod{12}$ for each $l$. In view of (\ref{eq:21}), (\ref{eq:22}), (\ref{eq:23}) and (\ref{eq:24}), we have
\begin{align*}
    \prod_{l=1}^r t_{1+}(p_l^\alpha)&\equiv((-1)^{y_1+\cdots+y_r})^\alpha\prod_{l=1}^r\left((\alpha+1)x_l^{12\alpha}+6\alpha(\alpha+1)x_l^{12\alpha-2}y_l^2\right)\pmod{3^5}\\
    &\equiv((-1)^{y_1+\cdots+y_r})^\alpha(\alpha+1)^r\\
    &\quad\times\left((x_1\cdots x_r)^{12\alpha}+6\alpha(x_1\cdots x_r)^{12\alpha-2}\sum_{\substack{i_1<\cdots<i_{r-1}\\j\neq i_s}}(x_{i_1}\cdots x_{i_{r-1}}y_j)^2\right)\pmod{3^5},\\
    \prod_{l=1}^r t_{2+}(p_l^\alpha)&\equiv((-1)^{y_1+\cdots+y_r})^\alpha\prod_{l=1}^r\left((\alpha+1)x_l^{12\alpha}+123\alpha(\alpha+1)(11\alpha+26)x_l^{12\alpha-2}y_l^2\right)\pmod{3^5}\\
    &\equiv((-1)^{y_1+\cdots+y_r})^\alpha(\alpha+1)^r\\
    &\quad\times\left((x_1\cdots x_r)^{12\alpha}+123\alpha(11\alpha+26)(x_1\cdots x_r)^{12\alpha-2}\sum_{\substack{i_1<\cdots<i_{r-1}\\j\neq i_s}}(x_{i_1}\cdots x_{i_{r-1}}y_j)^2\right)\pmod{3^5}
\end{align*}
Thus, going by the notations of (\ref{eq:29}), (\ref{eq:27}) and (\ref{eq:28}), we get
\begin{align*}
    &\prod_{\substack{p\mid(12n+13)\\p\equiv1(12)}}t_{1+}(p^\alpha)\equiv\prod_{\substack{p\mid(12n+13)\\p\equiv1(12)}}t_{2+}(p^\alpha)\pmod{3^5}\\
    \Longleftrightarrow&\prod_\alpha(\alpha+1)^{r_\alpha}(A_\alpha+6\alpha B_\alpha)\equiv\prod_\alpha(\alpha+1)^{r_\alpha}(A_\alpha+123\alpha(11\alpha+26)B_\alpha)\pmod{3^5},
\end{align*}
where the product is taken over all distinct $\alpha=ord_p(12n+13)$ with $p\equiv1\pmod{12}$. Since $B_\alpha\equiv0\pmod9$ for each $\alpha$, therefore expanding the products on both sides involving the terms $A_\alpha, B_\alpha$ and comparing we obtain
\begin{align*}
    &\prod_{\substack{p\mid(12n+13)\\p\equiv1(12)}}t_{1+}(p^\alpha)\equiv\prod_{\substack{p\mid(12n+13)\\p\equiv1(12)}}t_{2+}(p^\alpha)\pmod{3^5}\\
    \Longleftrightarrow&6\left(\prod_\alpha(\alpha+1)^{r_\alpha}\right)\left(\sum_\alpha\left(\alpha B_\alpha\prod_{\beta\neq\alpha}A_\beta\right)\right)\\
    &\qquad\qquad\equiv123\left(\prod_\alpha(\alpha+1)^{r_\alpha}\right)\left(\sum_\alpha\left(\alpha(11\alpha+26)B_\alpha\prod_{\beta\neq\alpha}A_\beta\right)\right)\pmod{3^5}\\
    \Longleftrightarrow&2\left(\prod_\alpha(\alpha+1)^{r_\alpha}\right)\left(\sum_\alpha\left(\alpha B_\alpha\prod_{\beta\neq\alpha}A_\beta\right)\right)\\
    &\qquad\qquad\equiv41\left(\prod_\alpha(\alpha+1)^{r_\alpha}\right)\left(\sum_\alpha\left(\alpha(11\alpha+26)B_\alpha\prod_{\beta\neq\alpha}A_\beta\right)\right)\pmod{3^4}
\end{align*}

Assume $k\equiv\varepsilon\pmod3$ in (\ref{eq:25}) where $\varepsilon=1$ or $2$. We have
\begin{align*}
    &b_{27}(n)\equiv0\pmod3\\
    \Longleftrightarrow&\prod_{\substack{p\mid(12n+13)\\p\equiv1(12)}}t_{1+}(p^\alpha)\prod_{\substack{p\mid(12n+13)\\p\equiv5(12)}}(-1)^{\alpha/2}p^{6\alpha}\equiv\prod_{\substack{p\mid(12n+13)\\p\equiv1(12)}}t_{2+}(p^\alpha)\left(81\varepsilon+\prod_{\substack{p\mid(12n+13)\\p\equiv5(12)}}(-1)^{\alpha/2}p^{6\alpha}\right)\pmod{3^5}\\
    \Longleftrightarrow&\prod_{\substack{p\mid(12n+13)\\p\equiv1(12)}}t_{1+}(p^\alpha)\prod_{\substack{p\mid(12n+13)\\p\equiv5(12)}}(-1)^{\alpha/2}p^{6\alpha}\equiv(81d+1)\prod_{\substack{p\mid(12n+13)\\p\equiv1(12)}}t_{2+}(p^\alpha)\prod_{\substack{p\mid(12n+13)\\p\equiv5(12)}}(-1)^{\alpha/2}p^{6\alpha}\pmod{3^5}\\
    \Longleftrightarrow&\prod_{\substack{p\mid(12n+13)\\p\equiv1(12)}}t_{1+}(p^\alpha)\equiv(81d+1)\prod_{\substack{p\mid(12n+13)\\p\equiv1(12)}}t_{2+}(p^\alpha)\pmod{3^5}\\
    \Longleftrightarrow&2\left(\prod_\alpha(\alpha+1)^{r_\alpha}\right)\left(\sum_\alpha\left(\alpha B_\alpha\prod_{\beta\neq\alpha}A_\beta\right)\right)\\
    &\qquad\qquad\equiv41(81d+1)\left(\prod_\alpha(\alpha+1)^{r_\alpha}\right)\left(\sum_\alpha\left(\alpha(11\alpha+26)B_\alpha\prod_{\beta\neq\alpha}A_\beta\right)\right)\pmod{3^4}
\end{align*}
where $d=\varepsilon d'$, $d'\in\mathbb Z$ is the inverse of $\prod_{\substack{p\mid(12n+13)\\p\equiv5(12)}}(-1)^{\alpha/2}p^{6\alpha}$ in $(\mathbb Z/3^5\mathbb Z)^*$.

\section*{\textbf{Acknowledgement}}

The author is thankful to Prof. Srilakshmi Krishnamoorthy for introducing this particular area of research. The author is also thankful to Prof. David Penniston for his invaluable initial guidance regarding the scope of research in this area.

\end{document}